\documentclass[11pt]{amsart}
\usepackage{amssymb}

\renewcommand{\bar}[1]{\overline{#1}}

\newcommand{\EE}{\mathbb{E}} 
\newcommand{\eps}{\varepsilon}

\newcommand{\Cald}{\Pi_P} 
\newcommand{\rest}[1]{\downharpoonright_{#1}} 

\newcommand{\df}[1]{\mathfrak{#1}}
\newcommand{\curly}[1]{\mathcal{#1}}

\newcommand{\lrpar}[1]{\left( #1 \right)}
\newcommand{\lrspar}[1]{\left[ #1 \right]}
\newcommand{\lrbrac}[1]{\lbrace #1 \rbrace}

\renewcommand{\bar}[1]{\overline{#1}}

\renewcommand{\tilde}[1]{\widetilde{#1}} 
\renewcommand{\hat}[1]{\widehat{#1}} 

\newcommand\inv{\operatorname{inv}}

\usepackage[all]{xy}

\newcommand\mhy{---}

\newcommand\boxb[1]{\square_b}

\newcommand\wt{\widetilde}
\newcommand\lt{\operatorname{lt}}

\numberwithin{equation}{section}
\newcommand\paperintro%
        {%
         }
\newcommand\paperbody%
        {%
         }

\newtheorem{non-theorem}{Non-Theorem}

\newtheorem{theorem}{Theorem}
\newtheorem{corollary}[theorem]{Corollary}
\newtheorem{lemma}[theorem]{Lemma}
\newtheorem{proposition}[theorem]{Proposition}
\numberwithin{theorem}{section}

\theoremstyle{remark}

\newcommand\Kc{\operatorname{K_{\text{c}}}}

\newcommand\Kto{\operatorname{K^1}}
\newcommand\Kte{\operatorname{K^0}}
\newcommand\Kcp[1]{\operatorname{K^{#1}_{\text{c}}}}
\newcommand\Kco{\Kcp1}
\newcommand\Kce{\Kcp0}

\newcommand\scat{\operatorname{sc}}

\newcommand\bo{\operatorname{b}}

\newcommand\cu{\operatorname{cu}}

\newcommand\coF{{}^{\mathcal{C}}\kern-2pt\Lambda}

\newcommand\cFTs{{}^{\Phi}\overline{T}\kern-1pt{}^*}

\newcommand\sus{\operatorname{sus}}

\newcommand\heta{\hat\eta}

\newcommand\GL{\operatorname{GL}}

\newcommand\tG{\tilde{G}}
\newcommand\tV{\tilde{V}}
\newcommand\tW{\tilde{W}}

\newcommand\ie{i\@.e\@. }

\newcommand\cF{\mathcal{F}}

\newcommand\cH{\mathcal{H}}
\newcommand\cI{\mathcal{I}}

\newcommand\cK{\mathcal{K}}

\newcommand\cN{\mathcal{N}}

\newcommand\NN{\mathbb N}

\newcommand\RR{\mathbb R}

\newcommand\ZZ{\mathbb Z}

\newcommand\bbC{\mathbb C}

\newcommand\bbE{\mathbb E}
\newcommand\bbF{\mathbb F}

\newcommand\bbR{\mathbb R}
\newcommand\bbS{\mathbb S}

\newcommand\cS{\mathcal S}

\newcommand\CI{{\mathcal{C}}^{\infty}}

\newcommand\cFNs{{}^{\Phi}\overline N\kern-1pt{}^*}

\newcommand\Dom{\operatorname{Dom}}

\newcommand\Id{\operatorname{Id}}

\newcommand\Ran{\operatorname{Ran}}

\newcommand\dCI{\dot{\mathcal{C}}^{\infty}}

\newcommand\pa{\partial}

\newcommand\Mfor{\text{ for }}

\newcommand\Mst{\text{ s.t. }}

\newcommand\Mwith{\text{ with }}

\begin{document}
\title[Fredholm realizations]
{Fredholm realizations of elliptic symbols\\ on manifolds with boundary}

\author{Pierre Albin}
\author{Richard Melrose}
\address{Department of Mathematics, Massachusetts Institute of Technology}
\email{pierre@math.mit.edu}
\email{rbm@math.mit.edu}
\thanks{The first author was partially supported by an NSF postdoctoral
  fellowship and the second author received 
  partial support under NSF grant DMS-0408993.}

\begin{abstract} We show that the existence of a Fredholm element of the
  zero calculus of pseudodifferential operators on a compact manifold with
  boundary with a given elliptic symbol is determined, up to stability,
  by the vanishing of the Atiyah-Bott obstruction. It follows that, up to
  small deformations and stability, the same symbols have Fredholm
  realizations in the zero calculus, in the scattering calculus and in the
  transmission calculus of Boutet de Monvel.
\end{abstract}
\maketitle

\tableofcontents

\paperintro
\section*{Introduction} 

Consider a compact manifold with boundary, $X,$ or more generally a
fibration 
\begin{equation}
\xymatrix{X\ar@{-}[r]&M\ar[d]\\
&B}
\label{fresmb0sc.7}\end{equation}
with typical fibre $X.$ To a symbol
\begin{equation*}
 	a\in \CI\lrpar{S^*\lrpar{M/B};\hom\lrpar{E,F}}
\end{equation*}
for vector bundles $E,$ $F$ over $M,$ one can associate, after a little
`preparation' various choices of families of pseudodifferential operators
with symbol $a.$ In particular, we may consider operators of `zero' type in
$\Psi^0_0\lrpar{M/B;E,F}$ or of `scattering' type in
$\Psi^0_{\mathrm{sc}}\lrpar{M/B;E,F},$ or of `transmission' type in Boutet
de Monvel's algebra $\Psi^0_{\mathrm{Tr}}\lrpar{M/B;E,F}.$ In these, and
other, cases we can seek the conditions on $a$ under which it has a
realization, in the corresponding sense, which is Fredholm on the natural
$L^2$-spaces. The primary condition is always symbolic ellipticity, that
$a^{-1} \in \CI\lrpar{S^*\lrpar{M/B};\hom\lrpar{F,E}}$ should exist. For
the three `quantization procedures' mentioned above such a Fredholm
realization of an elliptic symbol exists if and only if the K-class
associated to the symbol, in terms of absolute K-theory with compact
supports, lies in the image of the K-theory supported in the interior
\begin{equation*}
 	\Kc\lrpar{T^*\lrpar{\mathrm{int}\lrpar{M}/B}} \ni \alpha
	\mapsto i_*\lrpar{\alpha}=[a]\in\Kc\lrpar{T^*\lrpar{M/B}}.
\end{equation*}
This is established here for the zero calculus (as introduced in
\cite{Mazzeo:Hodge}, \cite{Mazzeo-Melrose}) by computing K-theory
invariants of the algebra. For the scattering calculus it is shown in
\cite{Melrose-Rochon} and for the transmission case already in essence by
Boutet de Monvel in \cite{Boutet} (see the much more recent work by Melo,
Nest, and Schrohe \cite{MNS} and Melo, Schick, and Schrohe \cite{MSS});
this in turn is an extension of the results of Atiyah and Bott
\cite{AtiyahBott}.

The elements of the Lie algebra, $\curly{V}\lrpar{X}$, of smooth vector
fields on a compact manifold with boundary act on $\CI\lrpar{X}$. If
$0<\nu\in\CI\lrpar{X;\Omega}$ is a smooth density and $d\lrpar{V} \in
\CI\lrpar{X}$ is the divergence of $V$ with respect to $\nu$ then
\begin{equation*}
 	\int_X \lrpar{Vu}\cdot w \nu
	=\int_X u\cdot\lrpar{-V+d\lrpar{V}w} \nu
	+\int_{\pa X} n_V u\cdot w \nu'
\end{equation*}
where $n_v\in\CI\lrpar{\pa X}$ is the normal component of $V$ at $\pa Z$
computed with respect to a density $0<\nu'\in\CI\lrpar{\pa X;\Omega}.$
Thus, in order to get an action of the associated enveloping algebra, one
needs either to add boundary conditions, as in the transmission algebra, or
restrict to elements $V \in \curly{V}_{\text{b}}\lrpar{X} \subset
\curly{V}\lrpar{X}$ of the Lie algebra of vector fields tangent to the
boundary for which $n_V=0.$

Associated to $\curly{V}_{\text{b}}$ is an algebra of pseudodifferential
operators $\Psi_{\text{b}}^{\ZZ}\lrpar{X}$ containing the enveloping algebra
$\mathrm{Diff}_{\text{b}}^{\NN}\lrpar{X} \supset
\curly{V}_{\text{b}}\lrpar{X}$ and acting continuously, with adjoints, on
$\CI\lrpar{X}.$ The same is true for various Lie subalgebras of
$\curly{V}_{\text{b}}\lrpar{X}$. Here we are most interested in the natural
ones:
\begin{equation}
\begin{split}
	\curly{V}_0\lrpar{X} &=
	\{ V \in \curly{V}_{\text{b}}\lrpar{X} : V=xW,
	W\in\curly{V}\lrpar{X} \} \\ 
	\curly{V}_{\mathrm{sc}}\lrpar{X} &=
	\{ V \in \curly{V}_{\text{b}}\lrpar{X} : V=xW,
	W\in\curly{V}_{\text{b}}\lrpar{X} \} \\ 
\end{split} 
\label{fresmb0sc.22}\end{equation}
fixed in terms of a boundary defining function $x\in\CI\lrpar{X}$ on which
they do not depend. There are other interesting algebras which depend on
some choice of, or the existence of, more structure. These include
\begin{equation*}\begin{split}
	\curly{V}_{\mathrm{cu}}\lrpar{X} &=
	\{ V \in \curly{V}_{\text{b}}\lrpar{X} : V\cdot x \in
	x^2\CI\lrpar{X} \} \\ 
	\curly{V}_{\operatorname{\psi-\mathrm{b}}}\lrpar{X} &=
	\{ V \in \curly{V}_{\text{b}}\lrpar{X} : V\cdot\CI_{\psi}\lrpar{X}
	\subset x\CI\lrpar{X} \} \\ 
	\curly{V}_{\operatorname{\psi-\mathrm{cu}}}\lrpar{X} &=
	\curly{V}_{\mathrm{cu}}\lrpar{X} \cap
	\curly{V}_{\operatorname{\psi-\mathrm{b}}}(X) \\
	\curly{V}_{\theta}\lrpar{X} &=
	\{ V \in \curly{V}_0\lrpar{X} : V\cdot\wt\theta \in x^2\CI\lrpar{X} \}
\end{split}\end{equation*}
where $\pa X \xrightarrow{\psi} Y$ is a fibration of the boundary, $x \in
\CI\lrpar{X}$ is a boundary defining function (on which for instance
$\curly{V}_{\mathrm{cu}}$ depends) and $\wt\theta \in
\CI\lrpar{X,\Lambda^1}$ is a real $1$-form with $\theta = i_{\pa
X}^*\wt\theta$ non-vanishing.

In all of these cases there is a well-defined algebra of pseudodifferential
operators, denoted $\Psi^{\ZZ}_S\lrpar{X}$, for the `structures'
$S=\mathrm{b},$ $0,$ $\mathrm{sc},$ $\cu,$
$\operatorname{\psi-\mathrm{b}},$ $\operatorname{\psi-\mathrm{cu}},$
$\theta$ which includes the corresponding enveloping algebra.

Since the operators of order zero act on $\CI\lrpar{X}$ as a $*$-algebra
and extend to bounded operators on $L^2\lrpar{X}$ the $C^*$-closure exists
\begin{equation*}
	\bar{\Psi^0_S} \subset \curly B \lrpar{L^2\lrpar{X}},
\end{equation*}
so the K-theory of this, and related algebras, is defined.
In all cases
\begin{equation*}
	\dot{\Psi}^{-\infty}\lrpar{X} \equiv \dot{C}^{\infty}\lrpar{
	  X^2;\Omega_R},
\end{equation*}
the space of smoothing operators with kernels vanishing to infinite order
at the boundary, is an ideal. It has closure $\df{K}$, the ideal of
compact operators and so leads to the sequence
\begin{equation*}
	\df{K}\longrightarrow \bar{\Psi_S^0}\lrpar{X}\longrightarrow 
	\curly{A}_S^0\lrpar{X} = \bar{\Psi_S^0}\lrpar{X}/\df{K}
\end{equation*}
and the associated $6$-term exact `index sequence' in K-theory.

It is also possible to construct a smooth version of K-theory,
without passing to the $C^*$-closure, giving however essentially the same result (see \S\ref{Cstar}). In
some cases this is a little delicate since these algebras may
\emph{not} admit the functional calculus, in particular this is the
case for the `regular singular' algebras corresponding to $S=\mathrm{b},$
$0,$ $\operatorname{\psi-\mathrm{b}},$ $\theta.$ Nevertheless we may define
a `geometric' replacement for $K_1\lrpar{\curly{A}_S^0}$ as in
\cite{Melrose-Rochon}, which we denote $\cK_S^0\lrpar{X}$ and more
generally $\cK^0_S\lrpar{M/B}$ in the case of a fibration \eqref{fresmb0sc.7}
and, as necessary, possessing fibrewise structure $S.$ This K-group is
identified as the stable homotopy classes, up to bundle isomorphism, of the
elliptic and Fredholm elements of the modules $\Psi^0_S\lrpar{M/B;E,F}$
acting between sections of vector bundles. The analytic index for families
then becomes a homomorphism
\begin{equation*}
	\cK^0_S\lrpar{M/B}\xrightarrow{\mathrm{ind}_a} K\lrpar{B}.
\end{equation*}
Again in all the cases listed above (and others) there is a `symbolically
trivial' ideal
\begin{equation*}
	\Psi^{-\infty}_S\lrpar{X} \subset \Psi^0_S\lrpar{X},
\end{equation*}
which contains non-compact operators, and a corresponding smooth K-group,
$\cK^0_{-\infty, S}\lrpar{M/B},$ corresponding to Fredholm operators of the
form $\Id + A,$ $A \in \Psi^{-\infty}_S\lrpar{X}.$ The short exact sequence
of algebras
\begin{equation*}
\Psi^{-\infty}_S\lrpar{X}\overset{I}\longrightarrow \Psi^0_S\lrpar{X}
\longrightarrow \CI\lrpar{{}^SS^*X}[[\rho]]
\end{equation*}
with image the star algebra for symbols in terms of a quantization map,
induces a 6-term exact sequence in the smooth K-theory
\begin{equation}
	\xymatrix{
	\cK^0_{-\infty,S}\lrpar{M/B} \ar[r]^I & 
	\cK_S^0\lrpar{M/B} \ar[r]^(0.4)\sigma & 
	\Kce\lrpar{T^*\lrpar{M/B}} \ar[d]^{\beta^0} \\
	\Kco\lrpar{T^*\lrpar{M/B}} \ar[u]^{\beta^1} &
	\cK_S^1\lrpar{M/B} \ar[l]^(0.4)\sigma & 
	\cK^1_{-\infty,S}\lrpar{M/B} \ar[l]^I }  
\label{fresmb0sc.5}\end{equation}
although it should be emphasized that we do \emph{not} construct a general
theory from which this follows.

In terms of K-theory, the algebras again determine two extreme cases
which we call the \emph{universal} (for symbols) case and the
\emph{geometric} case, with all others `intermediate'. The universal
cases are characterized by the fact that the 6-term exact sequence above
decouples into two short exact sequences
\begin{equation*}
	S \text{ universal} \iff
\cK^i_{-\infty,S}(M/B)=\operatorname{K}_{\operatorname{c}}^i(B\times\bbR),\
\beta^0=\beta^1=0.
\end{equation*}
\begin{theorem}[\cite{Melrose-Rochon}, see Corollary \ref{bUniversal}
    below] The cusp and b structures are universal.\end{theorem}

The second extreme case is characterized by
\begin{equation}
\begin{gathered}
\begin{aligned}
	S &\text{ is geometric} \iff \\
&\text{the $6$-term exact sequence \eqref{fresmb0sc.5}
	  reduces to the geometric sequence:}
\end{aligned}\\
	\xymatrix{
	\Kco\lrpar{T^*\lrpar{\pa M/B}} \ar[r] & 
	\Kce\lrpar{T^*\lrpar{M/B}, T^*_{\pa M}\lrpar{M/B}} \ar[r] & 
	\Kce\lrpar{T^*\lrpar{M/B}} \ar[d] \\
	\Kco\lrpar{T^*\lrpar{M/B}} \ar[u] &
	\Kco\lrpar{T^*\lrpar{M/B}, T^*_{\pa M}\lrpar{M/B}} \ar[l] & 
	\Kce\lrpar{T^*\lrpar{\pa M/B}}\ar[l] } 
\end{gathered}
\label{fresmb0sc.38}\end{equation}
in which the connecting homomorphisms are geometric boundary maps (and the
K-theory of manifolds with boundary is absolute).

\begin{theorem}[\cite{Melrose-Rochon} for sc, Proposition \ref{TopSixTerm}
    below for zero]
The structures $0$ and $\scat$ are geometric.\end{theorem}

This is the main technical result of this paper, it answers the question
posed by the title for the zero calculus in the form:

\begin{proposition}\label{answer} An elliptic symbol
$a\in\CI\lrpar{S^*\lrpar{M/B};\hom\lrpar{E,F}}$ has a realization as the
symbol of a Fredholm element of $\Psi^0_0\lrpar{M/B;E,F}$ if and only if
for some $N$ there exists 
\begin{multline*}
\wt{a}\in\CI\lrpar{\bar{T^*\pa M/B};\hom\lrpar{E\oplus\bbC^N,F\oplus\bbC^N}}
\Mwith\\
\wt{a}^{-1}\in
\CI\lrpar{\bar{T^*\pa M/B};\hom\lrpar{F\oplus\bbC^N,E\oplus\bbC^N}}
\label{fresmb.4}\end{multline*}
 such that
\begin{equation*}
 	a\oplus\Id_{\bbC^N} \rest{S^*\pa M/B}= \wt{a}\rest{S^*\pa M/B}.
\end{equation*}
That is, $a$ can be (stably) trivialized over the boundary.
\end{proposition}

In the {\em universal} case as opposed to the geometric cases there is no
restriction on the symbol other than invertibility (hence the name).

The `intermediate' case $S=\operatorname{\psi-\mathrm{cu}}$ for $\psi$ a
non-trivial fibration of the boundary (so both fibre and base 
having positive dimension) is discussed in \cite{Melrose-Rochon}. Using an
`adiabatic limit' construction it will be shown in a subsequent paper that 
there are natural isomorphisms
\begin{equation*}
 	\cK^i_{\operatorname{\psi-\mathrm{cu}}}\lrpar{M/B} \longleftrightarrow
 	\cK^i_{\operatorname{\psi-\mathrm{b}}}\lrpar{M/B} \longleftrightarrow
	KK^i_B\lrpar{C_{\psi}\lrpar{M}, C\lrpar{B}}
\end{equation*}
where the KK-functor (linear over $K(B))$ is applied to the
continuous functions on $M$ which are fibre constant for $\psi.$ The
boundary maps for the 6-term sequence then reduce to versions of the
Atiyah-Singer index map for the fibration of the boundary.

Before detailing the various sections of this paper, we mention some related literature.
The article \cite{Ditsche} was released shortly after the first version of this article and continues previous work of Lauter \cite{Lauter} regarding the $C^*$-algebras generated by operators in the zero calculus on a single manifold. In it, Ditsche computes the K-theory of the corresponding Calkin algebra. We show that this computation is equivalent to ours in \S\ref{Cstar} in the case of a single manifold, though for a family of manifolds one would have to consider a KK group as explained in \cite{Melrose-Rochon}.

The approach of Savin in \cite{Savin} deals with manifolds with non-trivial fibrations at the boundary. He shows that groups of elliptic operators, similar to those we study here, are isomorphic to the K-homology of the stratified manifold resulting from collapsing the fibration. In the geometric case treated here, Savin recovers the classical Atiyah-Bott obstruction as the obstruction to the existence of elliptic boundary conditions of Lopatinskii-Schapiro type.

Finally, we also mention \cite{Schulze1} where an extension of the transmission algebra is decribed so as to allow arbitrary elliptic operators to be quantized as Fredholm operators. This approach is extended in \cite{Schulze2} to remove the transmission condition on the principal symbol and in \cite{Schulze3} to allow for a non-trivial fibration at the boundary.

After briefly recalling the properties of the zero calculus in
\S\ref{Operators}, smooth K-theory is defined in \S\ref{Ktheory}, and in
\S\ref{BvsC} we discuss the relation between the b and cusp calculi, both
to review the universality, in the sense discussed above, of their K-theory
and to establish the closely related issue of the (weak) contractibility of
certain groups and semigroups. Then, in \S\ref{sec:FredRes}, this
contractibility is used to identify $\cK^0_{-\infty,0}\lrpar{M/B}.$ The six
term exact sequence and the identification of the group
$\cK^0_0\lrpar{M/B}$ is carried out in \S\ref{6term}. In the final
section the approach of Boutet de Monvel is reviewed and compared to the
scattering and zero calculi. The analogue of Proposition~\ref{answer} holds
for the transmission algebra, and corresponds to the obstruction to the
existence of local boundary conditions found by Atiyah and Bott, and we
show that the various lifts to $K^0\lrpar{T^*M/B,T^*_{\pa M}M/B}$ are
consistent, i\@.e\@., they yield the same index.

We are happy to acknowledge the helpful comments of the referee of this paper.

\paperbody

\section{Operators}\label{Operators}

\subsection{Zero calculus} \label{Background}

The zero calculus was defined and developed in \cite{Mazzeo:Hodge},
\cite{Mazzeo-Melrose}, \cite{Mazzeo:Edge}, \cite{Lauter}. In this section
we briefly recall the definitions and particularly those facts about the
{\em reduced normal operator} that we need below.

The space of $0$-vector fields is defined on any manifold with boundary by
\eqref{fresmb0sc.22}. In local coordinates, $x, z_i$ where $x$ is a
boundary defining function, and $z_i$ are coordinates along the boundary,
$\curly{V}_0$ is locally spanned by $x\pa_x,$ $x\pa_{z_i}.$ Zero
differential operators are elements of the enveloping algebra of
$\curly{V}_0$ and therefore are locally of the form
\begin{equation*}
	\sum_{j+|\alpha|\leq k}
	a_{j,\alpha}\lrpar{x,z}\lrpar{x\pa_x}^j\lrpar{x\pa_z}^\alpha
\end{equation*}
with arbitrary smooth coefficients. Such differential operators are
contained in the algebra of zero pseudodifferential operators. In the
references above the latter are identified as a space of conormal
distributions on a suitable modification of the double space $X^2$ (the
zero stretched double space $X^2_0).$ 

The space of zero pseudodifferential operators of order $k$ acting between
sections of the vector bundles $E$ and $F$ is denoted
$\Psi^k_0\lrpar{X;E,F}$.  There are two model operators. The first is the
principal symbol which corresponds to the short exact sequence
\begin{equation*}
\Psi^{k-1}_0\lrpar{X;E,F} \hookrightarrow
\Psi^k_0\lrpar{X;E,F} \overset{\sigma}\twoheadrightarrow
	\CI\lrpar{S^*X;\hom\lrpar{\pi^*E,\pi^*F}\otimes N_k}
\end{equation*}
The second is the normal operator which models the behavior at the boundary and
corresponds to a short exact sequence
\begin{equation}
x\Psi^k_0\lrpar{X;E,F}\hookrightarrow 
	\Psi^k_0\lrpar{X;E,F}\overset{N}\twoheadrightarrow
	\Psi^k_{\inv}\lrpar{\tG/\pa X;E,F}.
\label{fresmb0sc.6}\end{equation}
Here, $\tG\longrightarrow \pa X$ is a fibration of compact manifolds which
arises geometrically as the front face of the stretched double space; it is
discussed more extensively below. The fibres are closed quarter spheres
which combine to give a compactification of the bundle of (solvable) Lie
groups over $\pa X$ corresponding to the (canonical) $0$-structure. The
range space in 
\eqref{fresmb0sc.6}, which we characterize in more detail below, consists 
of distributional sections of a bundle over $\tG$ which are smooth away
from the identity section, where they have a classical conormal singularity,
and which vanish to infinite order at the boundaries. The product is just
the smooth family of convolution products, on 
the fibres, with `$\inv$' standing for this invariant,
i\@.e\@. convolution, product.

These two symbol homomorphisms, taken together, are analogous to the
principal symbol in the pseudodifferential calculus on a closed
manifold. Thus, for instance, a pseudodifferential operator acting on the
appropriate $L^2$-space is compact precisely when both symbols vanish; it
is Fredholm when both symbols are invertible. In the case of the
convolution algebra on the solvable group, invertibility is in the sense of
bounded operators on the natural Sobolev space and the inverse need not be
of the same form. In preparation for other formulations of the normal
operator we reinterpret the naturality of the zero structure and describe
alternative constructions of the bundle of groups and $\Psi_{\inv}.$

\subsection{Model zero structure}

Let $V$ be a real vector space of dimension $n$ and suppose that
$L^+\subset V'$ is a half-line in the dual, so 
\begin{equation}
\exists\ v'\in V'\Mst w'\in L^+ \iff w'=tv'\Mfor t>0.
\label{fresmb0sc.12}\end{equation}
With $V$ we may associate another vector space of the same dimension 
\begin{equation}
\tV=L\otimes V,\ L=\bbR\cdot L^+.
\label{fresmb0sc.13}\end{equation}
Thus, $\tV$ can be identified as the space of linear vector fields on $V$
(sums of products of linear functions and constant vector fields) where the
linear coefficients are in $L.$ As such it is a solvable Lie algebra,
more precisely if we choose a decomposition of $V=\bbR_x\times\bbR^{n-1}_y$
in which $L^+$ is spanned by $x,$ the dual of the first factor, then $\tV$
is spanned by 
\begin{equation}
x\pa_x,\ x\pa_{y_i}.
\label{fresmb0sc.14}\end{equation}
The span of the second class of vector fields is a subspace of dimension
$n-1$ which is well-defined as the center of the Lie algebra 
\begin{equation}
\tW\subset\tV.
\label{fresmb0sc.15}\end{equation}

The elements of $\tV$ are clearly complete as vector fields on $V$ so
choosing a point $p\in V^+,$ \ie where $L^+>0,$ then $V^+$ has a Lie group
structure with identity $p$ and Lie algebra $\tV$ and this is 
the unique connected and simply connected group (up to isomorphism) with
this Lie algebra. The centre $\tW$ integrates to a foliation of $V^+$ over
$(0,\infty)$ which has an induced multiplicative structure coming from the
short exact sequence of Lie groups 
\begin{equation}
L^\circ\longrightarrow G_{V,L^+}\longrightarrow (0,\infty),
\label{fresmb0sc.16}\end{equation}
\ie $G_{V,L^+}$ is a semidirect product. Choosing a `normal element' in
$\tV,$ complementary to $\tW$ and acting on it as the identity through the
Lie bracket, fixes an image of $(0,\infty)$ in $G_{V,L^+}$ and hence
identifies all the fibres with $\tW$ as a vector space (rather than just an
affine space). The compactification of $G_{V,L^+}$ to $\tG_{V,L^+},$ simply
the radial compactification of the closure in $V,$ gives a compact manifold
independent of choices, since linear transformations preserving the
boundary lift to be smooth on the compactification.

Applying these constructions to fibres of the tangent bundle to $X$ over
the boundary, with $L^+$ being the positive part of the conormal line,
gives a bundle of groups over $\pa X.$ Note that the front face of the zero
double space has interior isomorphic to the inward-point half of the
tangent bundle to the boundary, with a canonical interior `identity'
section. Namely, the interior is canonically the inward-pointing part of
the (projective) spherical normal bundle to the boundary of the
diagonal. This consists of vectors $(v_R,v_L),$ lifted from left and right,
modulo those tangent to the boundary of the diagonal - which are of the form
$(v',v')$ where $v'$ is tangent to the boundary - and modulo the $\bbR^+$
action. It follows that if $N$ is inward pointing then $(N,N)$ is a
well-defined interior point of the front face. For this reason the bundle
of group structures on the front face is well-defined. Note that to extend
this natural identity point to an origin for each of the tangential fibres
it is only necessary to choose a global inward-pointing vector field $N$
since then the points $(tN,N),$ $t>0,$ become origins of the abelian
fibres.

\subsection{Reduced normal operator, order $-\infty$}

To understand the structure of the normal operator better we use the
Fourier decomposition corresponding to the abelian fibres of the solvable Lie
group to `reduce' it to a parametrized family of operators on an interval.
In the case of operators of order $-\infty$ this was carried out explicitly
by Lauter, \cite{Lauter}, who characterized the range of the homomorphism
to this {\em reduced} normal operator. The normal operator in
\eqref{fresmb0sc.6} at each fibre takes values in
$\dot{\cS}([0,\infty)_s\times\bbR^{n-1}_Y),$ the space of Schwartz
functions on $\bbR^n$ with support in this half-space. The reduced normal
operator corresponds to Fourier transformation in $Y,$ followed by the
introduction of polar coordinates in the dual variable. The result
\begin{equation}
\cN(s,|\eta|,\heta)\in\cS([0,\infty)_s\times[0,\infty)\times\bbS^{n-2}_{\heta})
\label{fresmb0sc.8}\end{equation}
vanishes to infinite order at $s=0$ and may be interpreted as the kernel of an
operator on the $|\eta|$ half-line. Namely we may identify
$(0,\infty)\times(0,\infty)$ with the interior of
\begin{equation}
I^2_{\bo,0}=[I^2;\{(0,0)\}],\ I=[0,1]
\label{fresmb0sc.9}\end{equation}
by the map 
\begin{equation}
(s,r)\longrightarrow (sr,r)\in(0,\infty)^2\longrightarrow
  \lrpar{ \frac{sr}{1+sr},\frac{r}{1+r} } \in I^2.
\label{fresmb0sc.10}\end{equation}

In terms of the variables $s,$ $|\eta|,$ $\heta$ in \eqref{fresmb0sc.8},
$\cN(s,|\eta|,\heta)$ is an arbitrary smooth function which decays rapidly
with all derivatives as $s\to0,$ $s\to\infty$ and $|\eta|\to\infty,$ except
that there is special behaviour at $|\eta|=0$ where it is, of course, a
smooth function of $\eta$ not just smooth in $|\eta|$ and $\heta.$

Let $\pi$ be the projection $S^*\pa X\longrightarrow \pa X,$ let
$\lrpar{y,\eta}$ denote a point in $S^*\pa X,$ and let $\lrpar{\tau,\rho_b}
\in \lrspar{-1,1}\times\bbR^+$ be coordinates in 
the b, $\cu$-double space ($\rho_b=0$ defines the b-front face).

\begin{proposition}[\emph{cf} \cite{Lauter}, Prop. 4.4.1] \label{CharacRes} The
reduced normal operator defines an isomorphism of algebras
\begin{equation}
\Psi^{-\infty}_{\inv}(\tG/\pa X;E)\longleftrightarrow
\Ran^{-\infty}\left(\cN\right)\subset
\Psi_{\operatorname{b},\cu}^{-\infty,\cdot}(\cI/S^*\pa X;E)
\label{fresmb0sc.11}\end{equation}
where $\Psi_{\operatorname{b},\cu}^{-\infty,\cdot}(\cI/S^*\pa X;E)$ is
the space of b-pseudodifferential operators of order $-\infty,$ on the
compactified radial interval bundle, vanishing rapidly at the `infinite' (cusp)
end and $\Ran^{-\infty}\left(\cN\right)$ is the subspace with the additional
constraint that the Taylor series at the front face is of the form 
\begin{equation}
\sum\limits_{\alpha}I_\alpha (s)\eta ^\alpha.
\label{fresmb0sc.17}\end{equation}
\end{proposition}

\subsection{Reduced normal operator, general order}

In the general case of operators of integral order, the reduced normal
operator is defined in the same way but takes values in the
pseudodifferential operators on the unit interval which localize to an element of
the b-calculus near the left endpoint and to an element of the (weighted)
cusp calculus near the right endpoint with parameters in the cosphere
bundle of the boundary.
Indeed, if $A \in \Psi^k_0\lrpar{X}$, then
near the boundary the kernel of $A$ is given locally by
\begin{equation*}
	\curly{K}_A\lrpar{\tau,U,r,y}
	= \int_{\RR_\eta^{n-1}\times\RR_\xi} e^{i\lrpar{U\cdot\eta+\tau\xi}}a\lrpar{r,y,\xi,\eta} d\xi d\eta
\end{equation*}
for a symbol $a \in S^m_{cl}\lrpar{\RR_+\times\RR^{n-1}_y;\RR_\xi\times\RR^{n-1}_\eta}$.
We obtain the reduced normal operator by restricting to the zero front face ($r=0$) and Fourier transforming in the horizontal directions. 
The kernel of the reduced normal operator (as a half-density) is given by
\begin{equation}\label{ZeroKerN}
	\curly{K}_{\curly{N}\lrpar{A}}\lrpar{y,\eta;\tau,\rho_b}
	= \int_{\RR_\xi} e^{i\tau\xi}a\lrpar{0,y,\xi,\rho_b\eta} d\xi 
	\Bigl\lvert \frac{d\rho_b}{\rho_b}d\tau \Bigr\rvert^{\frac12}.
\end{equation}
This kernel defines, for fixed $\lrpar{y,\eta}$, a b,c-operator on the interval $\curly{I}=\lrspar{0,1}$.

This family of b,c-operators has, in turn, three model operators, the interior symbol corresponding to the 
conormal singularity at the diagonal, an indicial family at the b-end, and
a suspended family at the cusp end (see \cite{APSBook},
\cite{Mazzeo-Melrose1} respectively). 
The interior symbol is given at the point $\lrpar{u;\omega}\in S^*\curly{I}$ ($=\curly I\times\{\pm1\}$) by
\begin{equation}\label{bcsymbol}
	{}^{b,c}\sigma\lrpar{\curly{N}\lrpar{A}\lrpar{y,\eta}}\lrpar{u,\omega}
	= {}^0\sigma\lrpar{A}\lrpar{0,y;\omega,0} 
\end{equation}
as can be checked from \eqref{ZeroKerN} using the Taylor expansion of $a\lrpar{r,y,\xi,\eta}$ in $\eta$.
The kernel of the operator restricted to the b-face is
\begin{equation*}
	\int_{\RR_\xi} e^{i\tau\xi}a\lrpar{0,y,\xi,0} d\xi
\end{equation*}
and Mellin transform in $\tau$ produces the b-indicial family.

For the behavior near the cusp end, we first introduce the new coordinate
$\rho_{\cu}=\frac1\rho_{\operatorname{b}}$. The operator is then defined for
$\lrpar{\tau,\rho_{\cu}}\in\lrspar{-1,1}\times\lrspar{0,1}$, and the cusp front
face is obtained by blowing up the point $\lrpar{0,0}$. Lifting
$\curly{N}\lrpar{A}$ to the blown-up space (introducing
$T=\frac\tau{\rho_{\cu}}$) and restricting to a neighborhood near the new front
face and the diagonal, the kernel takes the form
\begin{equation*}\begin{split}
\curly{K}_{\curly{N}\lrpar{A}\lrpar{y,\eta}}\lrpar{s,T}
= &\int_{\RR_\xi} e^{i\rho_{\cu}T\xi}a\lrpar{0,y,\xi,\frac\eta{\rho_{\cu}}} d\xi 
\phantom{x}\rho_{\cu}\Bigl\lvert \frac{d\rho_{\cu}}{\rho_{\cu}^2}dT
\Bigr\rvert^{\frac12} \\
\xrightarrow{\hat\xi=\rho_{\cu}\xi}
&\int_{\RR_\xi} e^{iT\hat\xi}a\lrpar{0,y,\frac{\hat\xi}{\rho_{\cu}},
\frac\eta{\rho_{\cu}}} d\hat\xi 
\Bigl\lvert \frac{d\rho_{\cu}}{\rho_{\cu}^2}dT \Bigr\rvert^{\frac12}.
\end{split}\end{equation*}
Thus the restriction to the cusp front face has kernel
\begin{equation*}
\rho_{\cu}^{-k}\int_{\RR_\xi} e^{iT\hat\xi}
\phantom{x}{}^0\sigma\lrpar{A}\lrpar{0,y,\hat\xi,\eta} d\hat\xi
\end{equation*}
and the cusp suspended family is given by
\begin{equation*}
I_{\cu}\lrpar{\curly{N}\lrpar{A}\lrpar{y,\eta}}\lrpar{\xi}=
{}^0\sigma\lrpar{A}\lrpar{0,y,\xi,\eta}. 
\end{equation*}

In this way,
\begin{equation}
\cN:\Psi^k_0\lrpar{X;E,F}\longrightarrow 
	\rho_{\cu}^{-k}\Psi^k_{\operatorname{b},\cu}\lrpar{\cI/S^*\pa X;E,F}. 
\label{fresmb0sc.18}\end{equation}
These operators act on the fibres of the radial part of
the dual spaces to the abelian fibres in the compactified group. 
As in the smoothing case a key
feature of the reduced normal operator is that the b-indicial family only
depends on the parameters in $\pa X$ and not on the fibre in $S^*\pa X.$
Although we do not give an explicit characterization of the range of
\eqref{fresmb0sc.18} the following exact sequence suffices for our purposes.

\begin{proposition}\label{fresmb0sc.19} The (full) b-symbol map on the
  range $\Ran^{k}\left(\cN\right)$ of the reduced normal operator in
  \eqref{fresmb0sc.18} leads to a short exact sequence (with $k'=k+\frac n4$)
\begin{equation}
\rho_b^\infty\rho_{\cu}^{-k}
\Psi^{-\infty}_{\operatorname{b},\cu}(\cI/S^*\pa X;E,F)
\hookrightarrow
\Ran^{k}\left(\cN\right)
\twoheadrightarrow
I^{k'}_{\cS}(\cI_{\pa X})[[\eta]],
\label{fresmb0sc.20}\end{equation}
where the image space consists of the formal power series in $\eta$ with
the coefficient of $\eta ^\alpha$ being an arbitrary element of the space
  $I^{k'-|\alpha |}_{\cS}(\cI_{\pa X})$ of conormal distributions which are
  Schwartz at infinity and the null space is precisely the space of
  (b-)cusp operators of symbolic order $-\infty,$ vanishing to infinite
  order at the `b-end' and of singularity order $k$ at the cusp end (and in
  this sense arbitrary smooth sections of the bundle over $S^*\pa X).$
 \end{proposition}

\begin{proof} 
  From the discussion above, the full b-indicial operator is given by the
  Taylor expansion of the kernel at $\rho_b=0$: 
\begin{equation}\label{TaylorSeries}\begin{split}
	\cK_{\cN\lrpar{A}}\lrpar{y,\heta;\tau,\rho_b}
	&= \int_{\RR_\xi} e^{i\tau\xi}a\lrpar{0,y,\xi,\rho_b\heta} d\xi \\
	&\sim \sum_{k\geq 0} \sum_{|\alpha|=k} 
	\lrspar{ \int_{\RR_\xi} e^{i\tau\xi} D^\alpha a\lrpar{0,y,\xi,0}
  d\xi } \heta^\alpha \rho_b^k \\ 
	&\sim \sum_{|\alpha|\geq0} 
	\lrspar{ \int_{\RR_\xi} e^{i\tau\xi} D^\alpha a\lrpar{0,y,\xi,0}
  d\xi } \eta^\alpha 
\end{split}\end{equation}
  (with $\heta \in S^*\pa X$ and $\eta \in T^*\pa X$) hence has the form
  described in the theorem. 
  
  Conversely, any such formal power series can be asymptotically summed.
  Namely one can first get the
  symbol right, by taking Taylor series for the symbols, and then
  correcting the series for the smoothing terms. Thus the map is surjective.
  That is, given $\cK_\alpha \in I^{k'-|\alpha |}_{\cS}(\cI_{\pa X})$,
  we can find $a\lrpar{0,y,\xi,\eta}$ as in \eqref{TaylorSeries} such that
\begin{equation*}
	\int_{\RR_\xi} e^{i\tau\xi}a\lrpar{0,y,\xi,\eta} d\xi
	\sim
	\sum_{|\alpha| \geq 0} \cK_\alpha \eta^\alpha.
\end{equation*}
This condition only restricts the behavior of $a$ near
$\lrpar{0,y,\xi,\eta}$, hence the vanishing of the full b-symbol does not
restrict the possible cusp-symbols.  Finally, an element of the null space
has, by \eqref{bcsymbol}, $b,c$ symbol vanishing to infinite order, and so
is given by an (arbitrary) family in
\newline$\rho_b^\infty\rho_{\cu}^{-k}
\Psi^{-\infty}_{\operatorname{b},\cu}(\cI/S^*\pa
X;E,F).$
\end{proof}

\section{K-theory}\label{Ktheory}

\subsection{Smooth K-theory}

We define K-theory groups of zero operators as smooth versions of the
K-theory groups of the symbol algebra. Following \cite[Definition
 2]{Melrose-Rochon}, set
\begin{equation*}
A_0\lrpar{M;\bbE}=\{\lrpar{\sigma\lrpar{A}, \curly{N}\lrpar{A}}: A\in
\Psi^0_0\lrpar{M;\bbE} \phantom{x}\text{Fredholm on }L^2\} 
\end{equation*}
for a superbundle $\bbE=(E^+,E^-)$ and then define $\cK^0_0\lrpar{M}$ to
consist of the equivalence classes of elements of the union of the
$A_0\lrpar{M;\bbE}$ over $\bbE,$ where two elements are equivalent if there
is a finite chain consisting of one of the the following three equivalence
notions. First, 
\begin{align}
	A_{0}\lrpar{M;\bbE}\ni\lrpar{\sigma,N} &\sim \lrpar{\sigma',N'} 
		\in A_{0}\lrpar{M;\bbF} \label{Equiv1}\\
\intertext{if there are bundle isomorphisms $\bbE^{\pm}
\longrightarrow\bbF^{\pm}$ over $M$, that intertwine $\lrpar{\sigma,N}$ and
$\lrpar{\sigma',N'}.$ Secondly,}
 A_{0}\lrpar{M;\bbE}\ni\lrpar{\sigma,N} &\sim %
\lrpar{\wt\sigma,\wt{N}}\in A_{0}\lrpar{M;\bbE} \label{Equiv2}\\
	\intertext{if there exists a homotopy of Fredholm operators $A_t$ in 
		$\Psi_{0}^0\lrpar{M;\bbE},$ with
$\lrpar{\sigma,N} = \lrpar{\sigma\lrpar{A_0},\curly{N}\lrpar{A_0}}$ and
$\lrpar{\wt{\sigma},\wt{N}}=\lrpar{\sigma\lrpar{A_1},\curly{N}\lrpar{A_1}},$
	  and finally}
A_{0}\lrpar{M;\bbE}\ni\lrpar{\sigma,N} &\sim %
\lrpar{\sigma\oplus \Id_F, N \oplus \Id_F} \in A_{0}\lrpar{M;\bbE\oplus F}
\label{Equiv3}
\end{align}
where $\bbE\oplus F=(E^+\oplus F,E^-\oplus F).$ Similarly, we define
$\cK^1_{0}\lrpar{M}$ as consisting of the equivalence classes of elements
in the suspended space (corresponding to based loops)
\begin{equation*}
A_{0,\sus}\lrpar{M;E}
=\{s\in C^{\infty}\lrpar{\mathbb{R},A_{0}\lrpar{M;E}}:
s(t)-\Id_E\in\cS(\bbR;\Psi^0_{0}(X;E)\},
\end{equation*}
where $E=(E,E)$ as a superbundle and the equivalence condition is a finite
chain of relations \eqref{Equiv1}, \eqref{Equiv2}, \eqref{Equiv3} with bundle
transformations and homotopies required to be the identity to infinite
order at $\pm \infty$.

Similarly, we define $\cK_{-\infty,0}^0\lrpar{M}$ and
$\cK_{-\infty,0}^1\lrpar{M}$ as the groups of the corresponding equivalence
classes of elements in
\begin{equation*}
	A_{-\infty,0}\lrpar{M;E} =
	\{ \cN\lrpar{\Id + A} : A \in \Psi_{0}^{-\infty}\lrpar{M;E}, 
		\Id + A \phantom{x}\text{Fredholm on }L^2\},
\end{equation*}
and
\begin{equation*}
	A_{-\infty,0,\sus}\lrpar{M;E}
	= \{ s \in C^{\infty}\lrpar{\mathbb{R},A_{-\infty,0}\lrpar{M;E}} :
	s-\Id\in\cS(\bbR;\Psi^{-\infty}_0(X;E)\}
\end{equation*}
respectively.

\subsection{K-theory sequence}

From the compatibility of these two definitions and the standard definition
of (absolute) K-theory with compact supports it follows directly that there
is a sequence of maps 
\begin{equation}
\cK^i_{-\infty,0}(M/B)\longrightarrow \cK^i_{0}(M/B)\longrightarrow
\operatorname{K}^i_{\operatorname{c}}(T^*(M/B))
\label{fresmb0sc.24}\end{equation}
which is shown below to be exact.

\subsection{$C^*$ K-theory} \label{Cstar}

It is perhaps not immediately apparent that the definitions above lead to
groups, but this follows from the approximation properties discussed
below. 
In fact these also show that $\cK^i_0(X)$ is closely related to the K-theory of the
opposite parity of quotient of the $C^*$ closure of the algebra $\Psi^0_0(X)$ by the compact operators.
Although we will not use this below we discuss this connection in order to relate our results to those in the literature.

We start by fixing some notation, define
\begin{gather*}
	\cF_0(X;\bbE)
	= \lrbrac{ A \in \Psi^0_0(X;\bbE) : A \text{ is Fredholm on } L^2_0(X;\bbE) },
	\\
	\bar{ \cF_0(X;\bbE) }
	= \lrbrac{ A \in \bar\Psi^0_0(X;\bbE) : A \text{ is Fredholm on } L^2_0(X;\bbE) },
\end{gather*}
where $\bar\Psi^0_0(X;\bbE)$ is the $C^*$-closure of $\Psi^0_0(X;\bbE)$ as operators on $L^2_0(X;\bbE)$. Thus $A_0(X;\bbE)$ is just the symbol data of elements of $\cF_0(X;\bbE)$.
Also denote by $\df K$ the compact operators acting on $L^2_0(X)$.
The group $\cK_0^0(X)$ defined above is closely related to the group
\begin{gather*}
	K_1^{C^*}\lrpar{ \bar \Psi^0_0(X) / \df K } 
	= \lim_{\to} \mathrm{GL}_n( \bar \Psi^0_0(X) / \df K ) / \mathrm{GL}_n(\bar \Psi^0_0(X) / \df K )_0 \\
	= \lim_{\to} \bar{ \cF_0(X; \bbC^n) }/ \sim
\end{gather*}
where $\sim$ denotes homotopy.
Indeed this will follow immediately from the following lemma of Atiyah.

\begin{lemma}[\cite{AtiyahK}, Lemma A9]
Let $L\xrightarrow\pi M$ be a continuous linear map of Banach spaces with $\pi(L)$ dense in $M$ and let $U \subseteq M$ be open. Then for any compact $T$ the induced map
\begin{equation*}
	\lrspar{T, \pi^{-1}(U) } \to
	\lrspar{T, U } \text{ is bijective.}
\end{equation*}
\end{lemma}

{\emph Remark.} The proof of the lemma does not use that $L$ is complete, only that it is a normed vector space over a complete field.

\begin{proposition}
The group $K_1^{C^*}\lrpar{ \bar \Psi^0_0(X) / \df K }$ is isomorphic to  \newline 
$\displaystyle \lim_{\to} \cF_0(X; \bbC^n) / \sim$ and fits into the short exact sequence
\begin{equation*}
	0 \to K_1^{C^*}(C(X)) \to K_1^{C^*}\lrpar{ \bar \Psi^0_0(X) / \df K } \to \cK_0^0(X) \to 0.
\end{equation*}
\end{proposition}

\begin{proof}
It is clear that every equivalence class in $\bar{ \cF_0(X; \bbC^n) } / \sim$ has a representative from $\cF_0(X; \bbC^n)$ and Atiyah's lemma shows that elements of $\cF_0(X; \bbC^n)$ are homotopic in $\cF_0(X; \bbC^n)$ if and only if they are homotopic in $\bar{ \cF_0(X; \bbC^n) }$. Thus the groups $\cF_0(X; \bbC^n) / \sim$ and $\bar{ \cF_0(X; \bbC^n) }/ \sim$ are the same.

Exactness of the sequence follows from the fact that every non-compact manifold has a nowhere vanishing vector field (see e.g., \cite[Proposition 9]{MNS} for a proof).
In particular, for a Fredholm operator in $\Psi^0_0(X;\bbE)$ we always have $E^+ \cong E^-$. 
So by simultaneously stabilizing $E^+$ and $E^-$ we see that every class in $\cK_0^0(X)$ has a representative in $\cF_0(X;\bbC^n)$ and two representatives differ by a bundle isomorphism.
Finally, the existence of a non-vanishing section also implies that the map $K_1^{C^*}(C(X)) \to K_1^{C^*}(\bar \Psi/ \df K)$ is injective. Indeed, if $\Phi$ represents a class in $K_1^{C^*}(C(X))$ and $\Phi$ is homotopic to the identity through Fredholm operators in $\Psi^0_0(X;\bbC^n)$ then the principal symbol of any such homotopy evaluated at a non-vanishing section shows that $\Phi$ represents the identity in $K_1^{C^*}(C(X))$.
\end{proof}

A similar statement is true for the odd smooth K-theory group, $\cK^1_0(X)$. Our computations hold more generally for families of operators acting on a fibration $M \to B$ for which the analogue of the $C^*$ K-theory groups are the $K(B)$-linear $KK$ groups mentioned in the introduction.

\section{Relation between b and cusp calculi} \label{BvsC}

In this section, we briefly review from \cite{Surgery},
\cite{Melrose-Nistor} the relations between the b and cusp calculi. We
establish that the K-theory of the b-calculus is {\em universal} as
described in the introduction. For later use, we show that the group of
invertible elements in the $b, c$ calculus of order $-\infty$ vanishing to infinite order at the cusp end form
a contractible semigroup.

Recall that, on any manifold with boundary, $X,$ the Lie algebra of
b-vector fields, $\curly{V}_{\text{b}},$ consists of those vector fields
tangent to the boundary and in terms of local coordinates $x,$ $z_i$ ($x$ a
boundary defining function) near the boundary is spanned by $x\partial_x,$
$\partial_{z_i}.$ Similarly, for an admissible choice of $x$ for the cusp
structure, the Lie algebra of cusp vector fields, $\curly{V}_{\cu}(X)$ is
locally spanned by $x^2\partial_x,$ $\partial_{z_i}.$

The (canonical) b-structure on a manifold with boundary $X$ induces a
cusp-structure on $\wt{X}$, a manifold diffeomorphic to $X$ but with the
smooth structure enlarged to include the new boundary defining function 
\begin{equation*}
	r= \mathrm{ilg}\lrpar{x} = \frac1{\log \frac1x}.
\end{equation*}
Thus a function $f$ is smooth on $\wt{X}$ if it is smooth in the interior
and near the boundary is of the form
\begin{equation*}
	f\lrpar{r,z_i} = g\lrpar{\mathrm{ilg}\lrpar{x}, z_i}
\end{equation*}
for some function $g \in C^{\infty}\lrpar{[0,1)\times \pa X}.$ 
The identity map $\wt{X}\longrightarrow X$ is then smooth, corresponding to
the natural inclusion
\begin{equation*}
	C^{\infty}\lrpar{X}\longrightarrow C^{\infty}(\wt{X})
\end{equation*}
but in the opposite direction the identity map is not smooth. In this
sense, $\wt{X}$ has `more' smooth functions than $X$ (although it is
diffeomorphic to it, but not naturally so).  Both the interiors and the
boundaries of $X$ and $\wt{X}$ are canonically identified -- the manifolds
only differ in the manner in which the boundary is attached
(cf. \cite{Surgery}, \cite[$\S$2]{EpsteinMelroseMendoza}).

Now $\df{i}_{\mathrm{ilg}}^*\lrpar{x\pa_x} = r^2\pa_r$ and indeed,
$\df{i}_{\mathrm{ilg}}^*\curly{V}_{\text{b}}(X)$ spans
$\curly{V}_{\cu}(\wt{X})$ over $\CI(\wt{X}).$ The induced
relation between b and cusp differential operators extends to a close
relationship between the corresponding pseudodifferential operators. In fact, the
functor $X \to X_{\lt}$ (consisting of blowing-up all of the boundary faces
logarithmically and then `classically') takes the b-stretched double space
of $X$ to the cusp-stretched double space of $\wt{X},$ preserving
composition \cite[$\S$2.5-2.6, (4.4), $\S$4.6]{Surgery}. This leads to the
following result, the first part of which is
\cite[Prop. 26]{Melrose-Nistor0}. The second part is just the observation
that the lower order parts of the expansion of the kernel of an operator at
the front face can be smoothly deformed away without changing the normal
operator and preserving ellipticity.

\begin{proposition} \label{XN26} $ $
\begin{itemize}
\item [i)] 
The lift $(\beta^2_{lt})^*\Psi_{\text{b}}^s\lrpar{X}$ is a subcalculus of
$\Psi_{\cu}^s(\wt{X})$ with dense span over $C^{\infty}(\wt{X}).$ 
\item [ii)]
Any fully elliptic element of $\Psi^s_{\cu}(\wt{X})$ is homotopic through
fully elliptic elements to an element in
$\lrpar{\beta^2_{lt}}^*\Psi_{\text{b}}^s\lrpar{X}.$
\end{itemize}
\end{proposition}

\begin{corollary}\label{bUniversal}
The K-theory groups of the b and cusp calculi are naturally isomorphic,
\begin{equation*}
 	\cK^*_{-\infty,\mathrm b}\lrpar{M/B} \cong
 	\cK^*_{-\infty,\mathrm {cu}}\lrpar{M/B},\
 	\cK^*_{\mathrm b}\lrpar{M/B} \cong
 	\cK^*_{\mathrm{cu}}\lrpar{M/B};
\end{equation*}
since the cusp calculus is `universal' \cite{Melrose-Rochon}, so is the
b-calculus.
\end{corollary}

As in \cite{Melrose-Rochon0} we set
\begin{equation}\label{DefG}\begin{split}
\dot{G}^{-\infty}\lrpar{X;E}
	&= \{ \text{elements of } \Id  + \dot{\Psi}^{-\infty}\lrpar{X;E} \text{ invertible on $L^2$} \},\\
G^{-\infty}_{\mathcal S}\lrpar{X;E}
	&= \{ \text{elements of } \Id  + \Psi^{-\infty}_{\mathcal S}\lrpar{X;E} \text{ invertible on $L^2$} \},
\end{split}\end{equation}
and then the images of the normal operators on $G^{-\infty}_{\cu}\lrpar{X;E}$
and $G^{-\infty}_{\text{b}}\lrpar{X;E}$ are the corresponding subgroups on
which the index vanishes; they will be denoted 
\begin{equation*}
G_{\sus,\mathrm{ind}=0}^{-\infty}\lrpar{\pa X;E},\
G_{I,\mathrm{ind}=0}^{-\infty}\lrpar{\pa X;E}
\label{fresmb0sc.35}\end{equation*}
respectively. An element of the former can be smoothly deformed to an
element of the latter thanks to Proposition \ref{XN26}(ii).

In particular for an interval, we will also need the case of a manifold
with disconnected boundary, say
\begin{equation*}
	\pa X = \pa X_0 \sqcup \pa X_1
\end{equation*}
and the subgroup $G^{-\infty,\infty}_{\text{b}}\lrpar{X;E}$ in
$G^{-\infty}_{\text{b}}\lrpar{X;E}$ consisting of those operators which are
perturbations of the identity with
kernels vanishing to infinite order at the corner $\pa X_1 \times \pa X_1.$

\begin{corollary}\label{WeakContract}
The semigroups $G^{-\infty}_{\text{b}}\lrpar{X;E}$ and
$G^{-\infty,\infty}_{\text{b}}\lrpar{X;E}$ are weakly contractible.
\end{corollary}

\begin{proof}
The proof of \cite[Lemma 1]{Melrose-Rochon0}, that the sequence
\begin{equation*}
	\xymatrix{ 
	\dot{G}^{-\infty}\lrpar{X;E} \ar[r] 
	& G^{-\infty}_{\cu}\lrpar{X;E} \ar[r]^(.45){I_{\cu}} 
	& G_{\sus,\mathrm{ind}=0}^{-\infty}\lrpar{\pa X;E} }
\end{equation*}
is a Serre fibration, produces a lift, $h_t',$ of any family
\begin{equation*}
h_t: \curly{I}^k\longrightarrow 
G_{\sus,\mathrm{ind}=0}^{-\infty}\lrpar{\pa X;E}[[x]],\ t\in\lrspar{0,1}
\end{equation*}
to $G^{-\infty}_{\cu}\lrpar{X;E}$ with the convenient property that
$I_{\cu}\lrpar{h_t'}=h_t.$ 
This construction does not make use of the group structure of
$G^{-\infty}_{\cu}\lrpar{X;E}$ and so holds verbatim for the sequence
\begin{equation}\label{GbExact}
	\xymatrix{ 
	\dot{G}^{-\infty}\lrpar{X;E} \ar[r] 
	& G^{-\infty}_{\text{b}}\lrpar{X;E} \ar[r]^(.45){I_{\text{b}}} 
	& G_{I,\mathrm{ind}=0}^{-\infty}\lrpar{\pa X;E} }
\end{equation}
and for the restriction
\begin{equation}\label{GbbExact}
	\xymatrix{ 
	\dot{G}^{-\infty}\lrpar{X;E} \ar[r] 
	& G^{-\infty,-\infty}_{\text{b}}\lrpar{X;E} \ar[r]^(.45){I_{\text{b}}} 
	& G_{I,\mathrm{ind}=0}^{-\infty}\lrpar{\pa X_0;E} }
\end{equation}
when $\pa X = \pa X_0 \sqcup \pa X_1$.

Passage to $\wt{X}$ relates the two sequences by the commutative diagram
\begin{equation*}
\xymatrix{ 
\dot{G}^{-\infty}\lrpar{X;E} \ar[r] 
& G^{-\infty}_{\cu}\lrpar{X;E} \ar[r]^(.45){I_{\cu}} 
& G_{\sus,\mathrm{ind}=0}^{-\infty}\lrpar{\pa X;E} \\
\dot{G}^{-\infty}\lrpar{X;E} \ar[r] \ar[u]^{\df{i}^*} 
& G^{-\infty}_{\text{b}}\lrpar{X;E} \ar[r]^(.45){I_{\text{b}}} \ar[u]^{\df{i}^*}
& G_{I,\mathrm{ind}=0}^{-\infty}\lrpar{\pa X;E} \ar[u]^{\df{i}^*} },
\end{equation*}
and both the top and bottom rows induce long exact sequences of homotopy
groups. Since the left and right vertical arrows above are homotopy
equivalences, they induce isomorphisms of the corresponding homotopy
groups. The Five Lemma then shows that the middle vertical arrow is a weak
homotopy equivalence between $G^{-\infty}_{\text{b}}\lrpar{X;E}$ and the
contractible space $G^{-\infty}_{\cu}\lrpar{X;E}$. 
A similar argument shows the weak contractibility of
$G^{-\infty,-\infty}_{\text{b}}\lrpar{X;E}.$ 
\end{proof}

We remark that since
$G^{-\infty}_{\text{b}}\lrpar{X;E}$ is dominated by a CW-space, the
Whitehead lemma implies that $\df{i}^*$ is actually a homotopy equivalence
and hence that $G^{-\infty}_{\text{b}}\lrpar{X;E}$ is contractible \cite[Prop. 8.3]{Rochon}; we do
not use this in the sequel.

\section{Fredholm perturbations of the identity} \label{sec:FredRes}

In this section we identify the groups $\cK^*_{-\infty,0}\lrpar{M/B}$ with
the topological K-theory groups of the cotangent bundle of the boundary.

We quickly review the relevant definitions in terms of the groups
defined in \eqref{DefG}. Recall that for any manifold, $X$,
$K^1\lrpar{X}$ can be realized as stable homotopy classes of maps into
$\GL(N)$ or more directly (see \cite[Prop. 3]{Melrose-Rochon0}) in terms of
a classifying space
\begin{equation}
\Kto\lrpar{X} = \lim_{\to} \lrspar{X; \mathrm{GL}\lrpar{N} }
= \lrspar{X; \dot{G}^{-\infty}\lrpar{V;E} },
\label{K1def}\end{equation}
with $V$ any manifold with boundary. Correspondingly (see
\cite[Prop. 4]{Melrose-Rochon0}) 
\begin{equation}
\Kte\lrpar{X}=\lim_{\to} \lrspar{X; C^{\infty}\lrpar{
    \lrpar{\mathbb{S}^1, 1}; \lrpar{\mathrm{GL}\lrpar{N}, \Id}}} 
	= \lrspar{X; G^{-\infty}_{\sus}\lrpar{U;E}},
\label{K0def}\end{equation}
with $U$ any closed manifold of positive dimension. For a compact manifold
with boundary these definitions give the \emph{absolute} K-theory which we
denote (somewhat unconventionally) in the same way. The compactly supported
K-theory $\Kcp*(Y)$ of a non-compact manifold (possibly with boundary) is
defined in the same way, except that maps and homotopies are required to
reduce to the identity outside some compact subset of $Y.$ 

\begin{theorem}\label{0iso} For any fibration with fibres which are compact
  manifolds with boundary there is a canonical isomorphism
\begin{equation}
\cK^*_{-\infty,0}\lrpar{M/B}\longrightarrow \Kcp*\lrpar{T^*\pa(M/B)}.
\label{fresmb0sc.26}\end{equation}
\end{theorem}

\begin{proof} For simplicity of notation we first consider the case that
there is only a single compact manifold with boundary, $X.$ Given an
element $\kappa\in\cK^0_{-\infty,0}\lrpar{X},$ represented by $\Id + A$
with $A\in\Psi^{-\infty}_0\lrpar{X;E}$ we proceed to associate to it an
element $r(\kappa )\in\Kcp0\lrpar{T^*\pa X}.$ First, by adding the identity
on a complementary bundle, we may stabilize $E$ to $\bbC^M$ for some large
enough $M$ without changing $\kappa.$ In fact $\kappa$ is determined by the
normal operator of $\Id+A,$ and hence its reduced normal operator,
$\Id+\alpha.$

The condition that $\Id+A$ be Fredholm is equivalent to the invertibility
(on $L^2)$ of $\Id+\alpha$ as a family of operators on an interval bundle
over $S^*X.$ Proposition~\ref{CharacRes} characterizes the range of the
reduced normal operator in the calculus of b-pseudodifferential operators
of order $-\infty$ on the interval (and trivial at one end). Let us denote by
\begin{equation}
\cH^{-\infty}(X;E)\subset
\rho_{\cu}^\infty\Psi_{\operatorname{b},\cu}^{-\infty}(\cI/S^*\pa X;E)
\label{fresmb0sc.31}\end{equation}
the semigroup of those elements which are
invertible. Corollary~\ref{WeakContract} shows that the fibres of
$\cH^{-\infty}(X;E)$ are weakly contractible after stabilization. Thus, any
section of $\cH^{-\infty}(X;E)$ is homotopic to the identity section after
stabilization into $\cH^{-\infty}(X;\bbC^M)$ for sufficiently large $M.$ In
particular this is the case for $\alpha,$ so there exists a smooth curve 
\begin{multline}
\alpha':[0,1]\longrightarrow \cH^{-\infty}(X;\bbC^M) \\
\text{ such that }\alpha'(t)=\alpha \Mfor t<\epsilon ,\ \alpha'(t)=\Id\Mfor
t>1-\epsilon
\label{fresmb0sc.32}\end{multline}
for some $\epsilon >0.$ This is a family of b-pseudodifferential operators
on a bundle of intervals over $[0,1]\times S^*(\pa X)$ and of necessity the
family of b-indicial families $\gamma(t)=I_{\text{b}}(\alpha'(t))$ is also
invertible as a family of matrices. Thus associated with $\alpha '$ is a map 
\begin{equation}
\gamma \in\CI\lrpar{T^*\pa X, G^{-\infty}_{\sus}\lrpar{\{0\},\bbC^M}}.
\label{fresmb0sc.33}\end{equation}
Here we have used the fact from \eqref{fresmb0sc.11} that near $t=0$ the
b-indicial family is independent of the fibre variables of $S^*(\pa X)$
and near $t=1$ it is, by construction, equal to the identity. Thus, the
variable $t$ can be interpreted as a compactified radial variable for
$T^*(\pa X)$ giving \eqref{fresmb0sc.33}.

This defines the desired map
\begin{equation}
\cK^0_{-\infty,0}\lrpar{X}\ni\kappa \longmapsto [\gamma ]\in\Kcp0\lrpar{T^*\pa X}
\label{KInftyD}\end{equation}
which is to say that the K-class defined by $\gamma$ depends only on
$\kappa $ and not the intermediate choices. As in \cite[Prop 5.19]{Rochon},
the final class is independent of the trivialization of $E$ and, since any
two homotopies between sections of $\cH^{-\infty}(X;\bbC^M)$ are themselves
homotopic after stabilization, it is independent of the choice of
$\alpha'.$ Also, note that the product of two such elements $\Id + A,$ $\Id
+B$ is mapped to the product in $\Kcp0\lrpar{T^*\partial M}$ since we
can take the product of the homotopies.

Now, it remains to show that \eqref{KInftyD} is an
isomorphism. Surjectivity follows from the fact that every element in
$\Kce(T^*\pa X)$ may be represented by a family $\gamma(t)$ which is of the
form discussed above, so is an invertible family of indicial operators,
independent of the fibre variables near $t=0.$ From
Proposition~\ref{CharacRes} this may be quantized to an elliptic b,c
operator of the form $\Id+\alpha'(t)$ (with $\alpha '(t)$ trivial at the
cusp end) reducing to the identity in $t>1-\epsilon.$ The index bundle of
such a family (see the discussion below) is trivial, hence it may be
perturbed by a smoothing family (depending on all variables) to be
invertible and this may be chosen to be in the range of the reduced normal
operator for $t<\epsilon.$ The result at $t=0$ corresponds to a fully
elliptic operator $\Id+A$ in the zero calculus and hence a class $\kappa$
mapping to the given $\gamma.$
  
To see injectivity, consider a class $\kappa$ such that the construction
above leads to a family (so for an appropriate choice of initial
representative) which is homotopic to the identity, as a map on $T^*\pa X.$
Then after a small perturbation it can be arranged that this homotopy is a
family of the form of $\gamma,$ thus there is a 2-parameter family of
indicial operators $\gamma (t,s)$ on $[0,1]^2$ with $\gamma (t,0)$ the
indicial family and $\gamma (t,s)$ independent of the fibre variables for
$s<\epsilon$ and reducing to the identity if either $s>1-\epsilon$ or
$t>1-\epsilon.$ Again from Proposition~\ref{CharacRes} there is no
obstruction to `quantizing' this to a family of reduced normal operators,
for $s<\epsilon$ and to a family of b,c operators with more general
dependence on $S^*\pa X$ in $s>\epsilon.$ Such a family $\alpha '(t,s)$ may
be chosen to be the identity where $\gamma$ is the identity and to reduce to
the original homotopy at $s=0.$ Now, as an elliptic family, it has a
virtual index bundle over $[0,1]^2\times S^*\pa X.$ As a bundle this may be
realized as the difference of the null and conull bundles for $\alpha
'(t,s)(\Id-\Pi)$ where $\Pi$ is a finite rank self-adjoint smoothing
operator with range in $\dCI(\cI)$ and sufficiently large rank. It follows
that the index bundle is trivial and that the family may be made invertible
by a similar perturbation. At $t=0$ this gives a homotopy within the range
of the reduced normal operator to the identity. Thus the class $\kappa$ is
trivial.

Although this proof is written out only in the case of a single manifold it
is the dependence on the fibre variables in $S^*\pa X$ which is crucial and
it carries over with only notational changes to the general setting of a
fibration \eqref{fresmb0sc.7}. The argument for the odd K-groups involves
only the addition of a parameter.
\end{proof}

\section{Six term exact sequence}\label{6term}

We start out with a theorem very similar to \cite[Lemma
2.2]{Melrose-Rochon} and to the Atiyah-Bott analysis of the index problem
for boundary value problems \cite{AtiyahBott}. We will exploit the
similarity to \cite{Melrose-Rochon} to identify a six-term exact sequence
for the zero calculus with that of the scattering calculus.

\begin{theorem}\label{HtpyAway} For any fibration \eqref{fresmb0sc.7} (with
fibres modeled on a compact manifold with boundary) the smooth K-theory
$\cK^i_{0}(M/B)$ is represented by the Fredholm zero operators which are equal
to the identity near the boundary and hence there is a natural
isomorphism
\begin{equation}
\cK^0_0(M/B)\longrightarrow \Kcp*(T^*(\overset{\circ}{M}/B)).
\label{fresmb0sc.36}\end{equation}
\end{theorem}

\begin{proof} Consider an element $A\in\Psi^0_{0}(M/B;\bbE)$ which is fully
elliptic and hence represents an element of $\cK^0_0(M/B).$ The reduced
normal family of $A$ is an invertible (on $L^2)$ element of the b,c
calculus on an interval bundle over $S^*(\pa M/B).$ The symbol of $A$, $\sigma(A)$, is a
bundle isomorphism between $E^+$ and $E^-$ over the zero 
cosphere bundle of the fibers of $M/B$.
At the boundary, the inward normal fixes a section of this bundle so
$\sigma (A)$ at this section may be used to identify $E^+$ and $E^-$ at,
and hence near, the boundary. Thus, $\sigma (A)$ may be taken to be the
identity at this section. 

The indicial family at the b-front face, being invertible, then gives a
section of the isomorphism bundle of the lift of $E$ to the whole of the
compactified normal bundle to the boundary. By construction this section is
the identity at the inward-pointing end and is consistent with the symbol at
the outward-pointing end. So we may deform (within the b,c-calculus, 
not {\em a priori} within the range of the reduced normal operator) the b-indicial operator
keeping the consistency with an elliptic symbol at the inward-pointing end
until the b-indicial operator is the identity and the symbol is the
identity near both inward- and outward-pointing normal sections. 
Now, from Proposition~\ref{fresmb0sc.19} this homotopy may realized in the image of
the reduced normal operators of the zero calculus, so there is an elliptic curve
$A_t\in\Psi^0_0(M/B;\bbE)$ with $A_0=A$ and such that $A_1$ has \emph{full}
b-indicial operator equal to the identity. Initially this curve need not be
fully elliptic, but the reduced normal family is fully elliptic. Again by
standard index theory this can be modified to an invertible family, without
change near $t=0,$ by adding a term of order $-\infty$ with support in the
interior of the interval.

So, after this initial homotopy we conclude that $\cK^0_0(M/B)$ is
represented by fully elliptic families $A\in\Psi^0_0(M/B,\bbE)$ where $E^+$
and $E^-$ are identified near the boundary and with symbol which is the
identity near both normal directions and with reduced normal family having
full b-indicial family equal to the identity. Thus the reduced normal
family is in fact a section over $S^*(\pa(M/B))$ of the bundle of groups
with fibre the group 
\begin{equation}
G^{\cdot,-\infty}_{\cu}(\cI;E)
\label{fresmb0sc.37}\end{equation}
that is, trivial near the b-end. As explained in \S\ref{BvsC}, after stabilization this group is (weakly) contractible. Thus in fact
the reduced normal section may be (after stabilization) contracted to the
identity section. As above, this family may be lifted to a fully elliptic
homotopy of $A$ to a family which is the identity near the boundary and in
particular has reduced normal family equal to the identity.

Now the isomorphism \eqref{fresmb0sc.36} is just the usual Atiyah-Singer
isomorphism for the symbol.
\end{proof}

\begin{corollary}
With the smooth K-theory of the zero calculus identified with the compactly
supported K-theory of the fibre cotangent bundle of the interior of $M,$
the analytic index factors through the Atiyah-Singer index map for the
double $2M=M\cup M^-$
\begin{equation*}
\xymatrix@C=-15pt{\cK_0^0\lrpar{M/B} \ar[rr]^{\mathrm{ind}_a} \ar[d]^{\cong} & &
K\lrpar{B} \\
	\Kc\lrpar{T^*\lrpar{M/B},T^*_{\pa M}\lrpar{M/B}} \ar[dr]^{\cong} &
 &
	\Kc\lrpar{T^*\lrpar{2M/B}} \ar[u]^{\mathrm{ind}_{AS}} \\
&\Kc\lrpar{T^*\lrpar{2M/B},T^*\lrpar{M^-/B}} \ar[ur]^(.6){\sigma}}
\end{equation*}
\end{corollary}

\begin{proposition}\label{TopSixTerm}
The smooth K-theory groups of the zero calculus lead to the inner six term
exact sequence
\begin{equation}\label{Alg6T}
	\xymatrix{
K^0_c\lrpar{T^*\pa M/B} \ar[r]\ar@{<-}@<4pt>`l[dd] `d[dd] & 
K^0_c\lrpar{T^*M/B, T^*_{\pa M}M/B} \ar[rd] & 
\\
\cK^0_{-\infty,0}\lrpar{M/B} \ar[r]\ar[u] & 
\cK^0_0\lrpar{M/B} \ar[r]\ar[u] &
K^0_c\lrpar{T^*M/B} \ar[d] \\
 K^1_c\lrpar{T^*M/B} \ar[u] &
\cK^1_0\lrpar{M/B} \ar[l]\ar[d] &
\cK^1_{-\infty,0}\lrpar{M/B} \ar[l]\ar[d]\\
&
K^1_c\lrpar{T^*M/B, T^*_{\pa M}M/B} \ar[lu] & 
K^1_c\lrpar{T^*\pa M/B} \ar[l]\ar@{<-}@<4pt>`r[uu] `u[uu]}
\end{equation}
which is naturally isomorphic to the outer, geometric, sequence and thus
the zero calculus is `geometric', in the sense of the introduction.
\end{proposition}

\begin{proof}
The exactness of the inner diagram follows from the known exactness of the
outer diagram and the commutativity of the whole, so it only remains to
show the latter.

We can fit the maps from Theorems \ref{0iso} and \ref{HtpyAway} into the diagram
\begin{equation*}
\xymatrix{
\cK^0_{-\infty,0}\lrpar{M/B} \ar[r] \ar[d]^{\cong} &  
\cK^0_0\lrpar{M/B} \ar[r] \ar[d]^{\cong}  & 
K^0_c\lrpar{T^*M/B}  \ar@{<->}[d]^= \\
 K^0_c\lrpar{T^*\pa M/B} \ar[r] & 
K^0_c\lrpar{T^*M/B, T^*_{\pa M}M/B} \ar[r] & 
K^0_c\lrpar{T^*M/B}   }
\end{equation*}
in which the right square clearly commutes.
For the commutativity of the left hand square, note that the first map along the bottom is given by
\begin{equation*}
 \xymatrix{
K^0_c\lrpar{T^*\pa M/B} \ar[rr] \ar[rd]^{\cong} & &  
K^0_c\lrpar{T^*M/B, T^*_{\pa M}M/B} \\
& K^1_c\lrpar{T^*_{\pa M}M/B} \ar[ru]^{\alpha} & }
\end{equation*}
where $\alpha$ is defined by interpreting $T^*_{\pa M}M$ as a collar neighborhood of $T^*\pa M$ in $T^*M$,
and the same collar neighborhood idea is used in the identification of 
$\cK^0_0\lrpar{M/B}$ and $K^0_c\lrpar{T^*M/B, T^*_{\pa M}M/B}$.
Indeed, given a representative $A$ of a class in $\cK^0_0\lrpar{M/B}$, we found a homotopy of the normal operators starting at $\curly N\lrpar{A}$ and ending at the identity.
We can think of this homotopy as happening along a collar neighborhood by identifying it with $T^*\pa X\times \lrspar{0,1}_\varepsilon$ and using $\varepsilon$ as the homotopy variable.
Finally, note that the identification of $\cK^0_{-\infty,0}\lrpar{M/B}$ with $K^0_c\lrpar{T^*\pa M/B}$ was also by means of a homotopy
(this time thinking of $T^*\pa M/B$ as $S^*\pa M/B \times \bbR^+$ and running the homotopy in $\bbR^+$)
 and so we can see that the left square also commutes.

The same arguments show the commutativity of the corresponding squares
involving the odd K-theory groups. The two remaining squares, involving
$\Kc^*\lrpar{T^*M/B} \to \cK^{*+1}_{-\infty,0}\lrpar{M/B}$, commute by the
definition of this map.  \end{proof}

We point out that this proposition answers one of the questions from the introduction. Namely, given an elliptic symbol when is there a Fredholm family of operators in the zero calculus with that symbol? Since this is the same as asking if the symbol is in the image of the map
\begin{equation*}
	\cK_0^0(M/B) \to K_c^0(T^*M/B)
\end{equation*}
in the diagram above, exactness of the diagram at $K^0_c(T^*M/B)$ answers this question completely.

\begin{corollary}
A family of operators in $\Psi^0_0\lrpar{M/B;\EE}$ can be perturbed by a
family of operators in $\Psi^{-\infty}_0\lrpar{M/B;\EE}$ to be Fredholm if
and only if its interior symbol is in
\begin{equation*}
	\ker\lrpar{ K^0\lrpar{T^*M/B} \to K^1\lrpar{T^*_{\pa M}M/B}},
\end{equation*}
\ie, if its Atiyah-Bott obstruction vanishes. In which case, after
stabilization, there is a perturbation that makes it invertible.
\end{corollary}

\section{Relation to boundary value problems} 

Finally we relate the various approaches to quantization mentioned in the
introduction: the scattering calculus, the zero calculus, elliptic boundary
value problems as discussed by Atiyah and Bott, and the transmission
algebra of Boutet de Monvel. We will confine ourselves to relating the
approaches to the index problem, as a complete study of the relations
between these various calculi is more involved.

\subsection{Atiyah-Patodi-Singer versus local boundary conditions}\label{APSvsLS}

We start by recalling the general context of elliptic boundary value
problems, initially for a single differential operator of first order. See
for instance the exposition for the case of a single Dirac operator in
\cite{BoosW}.

Given an elliptic first order operator $P:\CI\lrpar{X;E}\to\CI\lrpar{X;F}$
on a compact manifold $X$ with non-trivial boundary, its kernel is infinite
dimensional and obtaining a Fredholm operator requires a restriction
reducing this to a finite dimensional space. The Calder\'on projector of
$P$, $\Pi_P$, is the orthogonal projection onto the space of Cauchy data 
of $P$,
\begin{equation*}
CD\lrpar{P}= \{ s\downharpoonright_{\pa X} : s \in \CI\lrpar{X;E}, Ps=0\}.
\end{equation*}
An elliptic boundary condition is a pseudodifferential operator of order
zero taking values in some auxiliary bundle $G,$ 
\begin{equation*}
\CI\lrpar{\pa X;E\downharpoonright_{\pa X}} \xrightarrow{R} \CI\lrpar{\pa X;G},
\end{equation*}
satisfying two conditions:
\begin{itemize}
\item [i)] $\displaystyle H^{(s)}\lrpar{\pa X;E\downharpoonright_{\pa X}} \xrightarrow{R^{(s)}} H^{(s)}\lrpar{\pa X;G}$ has closed range for every $s$
 \item [ii)] $\displaystyle \Ran\lrpar{\sigma\lrpar{R}} 
 	= \Ran\lrpar{\sigma\lrpar{R}\sigma\lrpar{\Cald}}
	\cong \Ran\lrpar{\sigma\lrpar{\Cald}}.$
\end{itemize}
Given an elliptic boundary condition the corresponding `realization' of
$P$, $P_R$, given by restricting $P$ to
\begin{equation*}
\Dom\lrpar{P_R}= \{s\in \Dom\lrpar{P}:R\lrpar{s\rest{\pa X}}=0  \}
\end{equation*}
is a Fredholm operator \cite[Thm 20.8, Prop. 18.16]{BoosW}.

Two types of elliptic boundary conditions stand out.  An elliptic condition
is of {\em generalized Atiyah-Patodi-Singer} or APS type if $R$ is an
orthogonal projection with the same principal symbol as $\Cald$. APS type
boundary conditions always exist, e.g., $R=\Cald.$

An elliptic boundary condition is {\em local} if $\sigma (R)$ is
surjective. In this case, condition $(i)$ above is automatic \cite[Lemma
20.9]{BoosW}. The index of a local elliptic boundary problem only depends
on the principal symbols of $P$ and $R$. This is not true for general,
e\@.g\@. APS type, elliptic boundary conditions.

In \cite{MelrosePiazza} and \cite{LoyaMelrose}, \cite{Loya}, the relation
between elliptic boundary conditions of APS type and Fredholm
b-pseudodifferential operators was analysed, leading to

\begin{theorem}[\cite{MelrosePiazza}, Lemma 5; \cite{Loya}, Theorem 1.1] 
Let $\eth_T$ be a generalized
Dirac-type operator on a manifold with boundary with boundary conditions of
APS type, then there is a b-operator of order $-\infty,$ $\curly{T},$ such that
$\eth + \curly{T} \in \Psi_{\text{b}}^1\lrpar{M;E}$ is a Fredholm operator
with the same index as $\eth_T.$
\end{theorem}

\subsection{Atiyah-Bott obstruction}\label{sec:AB}

For a differential operator, $P$, of general degree, $k$, the Cauchy data
involves the $k-1$ normal jet bundle at the boundary.  It is standard
(e.g., \cite[Chapter X]{Horm}, \cite[$\S$II.6]{BoosB}) to study this space
by passing to a family of ordinary differential operators formed from the
principal symbol. This family,
$P_{y,\heta},$ acts on half-lines and is parameterized by the cosphere
bundle of the boundary. The ellipticity of $P$ gives rise to a
decomposition of the null space of the $P_{y,\heta}$
\begin{equation*}
 	\curly M_{\lrpar{y,\heta}} = 
 	\curly M_{\lrpar{y,\heta}}^+ \oplus 
 	\curly M_{\lrpar{y,\heta}}^-
\end{equation*}
where $\curly M_{\lrpar{y,\eta}}^\pm$ consists only of exponential
polynomials involving $e^{i\lambda x}$ with
$\pm\mathrm{Im}\lrpar{\lambda}>0$ respectively. A local boundary condition
defines a bundle map 
\begin{equation*}
\beta^+:\curly M^+ \to \pi^*\lrpar{G}.
\label{fresmb0sc.39}\end{equation*}
and the Lopatinskii-Schapiro condition for ellipticity is that this map be
an isomorphism.

Atiyah and Bott \cite{Atiyah}, \cite{AtiyahBott} showed that the existence
of elliptic boundary conditions, in the sense of Lopatinskii-Schapiro, for
a given differential operator, $P,$ is topologically obstructed in that 
the symbol of $P$ is then in the kernel of the map
\begin{equation*}
 	\Kc\lrpar{T^*X} \xrightarrow{r_*} \Kc\lrpar{T^*_{\pa X}X}.
\end{equation*}
The isomorphism $\beta^+$ then determines a lifting of $[\sigma
  (P)]\in\Kc\lrpar{T^*X}$ to $[\sigma (P),\sigma
  (R)]\in K\lrpar{T^*X,T^*_{\pa X}X}.$  

\subsection{Boutet de Monvel transmission algebra}

Boutet de Monvel \cite{Boutet} constructed an algebra of pseudodifferential
which includes the generalized inverses of Lopatinskii-Schapiro elliptic
boundary problems. A general element of the algebra acts between bundles $E_1,$
$E_2$ over $X$ and bundles $F_1,$ $F_2$ over $\pa X$ and is given by 
\begin{equation}\label{BoutetMat}
\curly A= \begin{pmatrix} P_+ + G &  & K \\ &  & \\T &  & Q \end{pmatrix}:
\begin{matrix} \CI\lrpar{X;E_1} \\ \oplus \\
\CI\lrpar{\pa X;F_1} \end{matrix}\longrightarrow 
\begin{matrix} \CI\lrpar{X;E_2} \\ \oplus \\ \CI\lrpar{\pa X;F_2} \end{matrix}
\end{equation}
where $P$ is a pseudodifferential operator satisfying the transmission condition,
$P_+$ is its action on $\CI(X;E_1)$ by extension-as-zero and restriction,
$G$ is a `singular Green operator', $K$ and $T$ are potential and trace
operators, and $Q$ is a pseudodifferential operator on the boundary.
A differential operator $P$ with local boundary condition can be realized
as the operator 
\begin{equation}\label{BoutetLS}
 	\begin{pmatrix} P  \\ R  \end{pmatrix}
\end{equation}
in the Boutet de Monvel algebra.

We will say that $\curly A$ is elliptic if $P$ is elliptic and {\em fully
elliptic} if $\curly A$ defines a Fredholm operator on the natural Sobolev
spaces. The Lopatinski-Schapiro conditions on $R$ are equivalent to asking
that \eqref{BoutetLS} be fully elliptic. With each operator $\curly A$ as
in \eqref{BoutetMat} is associated a `boundary symbol' which is a family of
operators of Wiener-Hopf type. Ellipticity of $P$ implies that this family
is Fredholm and its index bundle is isomorphic to $\curly M^+$
\cite[pg. 35]{Boutet}. The operator $\curly A$ is fully elliptic if and
only if it is elliptic and its boundary symbol is invertible.

By \cite[Theorem 5.14]{Boutet}, an elliptic pseudodifferential operator $P$
satisfying the transmission condition has a realization as a Fredholm
operator $\curly A$ if and only if the Atiyah-Bott obstruction
vanishes. Boutet de Monvel showed that any such fully elliptic operator is
homotopic, through fully elliptic operators, to an element of the form
\begin{equation*}
\begin{pmatrix} P'_+ & 0 \\ 0 & Q' \end{pmatrix},
\end{equation*}
where, moreover, $P'_+$ is equal to the identity in a neighborhood of the
identity.

Melo, Schick, and Schrohe \cite{MSS}, following work of Melo, Nest, and
Schrohe \cite{MNS} showed that the K-theory of the $C^*$-closure of the
symbol algebra in the Boutet de Monvel calculus, $\cK_*$, gives a short
exact sequence
\begin{equation*}
	0 \to K_*\lrpar{C\lrpar{X}} \to \cK_* \to K_{*+1}\lrpar{T^*X^o} \to 0
\end{equation*}
and the analytic index map factors through $K_{*+1}\lrpar{T^*X^o}$.

Boutet de Monvel \cite[pg. 33]{Boutet} explicitly mentions an extension to
families of operators, proving that there exists an elliptic Green system
$\curly A_b$ (depending continuously on $b\in B$) associated to $P_b$ if
and only if $j\lrpar{P^+_b}\in K\lrpar{S^*\pa X \times B}$ is the pull-back
of a virtual bundle on $\pa X \times B.$

\subsection{Relationship to scattering and zero calculi}

The similarities between the three calculi mentioned in the introduction is
striking given the disparate motivations behind them.  Whereas the Boutet
de Monvel transmission algebra was constructed precisely to deal with
boundary value problems, the zero and scattering calculus serve to model
certain asymptotically regular, non-compact manifolds which admit a smooth
compactification to manifolds with boundary. The zero calculus is used to
model asymptotically hyperbolic geometries and the scattering,
asymptotically Euclidean. Nevertheless, the relation between the fully
elliptic elements of the calculi is very close.

\begin{theorem}
For an element $p \in \Kc\lrpar{T^*X}$ the following are equivalent:
\begin{enumerate}
 \item There is a fully elliptic element of the transmission algebra 
\begin{equation*}
 	\curly{A} = \begin{pmatrix} P_++G &  K \\ T  & Q \end{pmatrix}
	 \in \Psi^0_{\mathrm{tn}}\lrpar{X;\bbE,\bbF}
	\text{such that $[\sigma\lrpar{P_+}] = p.$}
\end{equation*}
  \item There is a fully elliptic element of the zero calculus 
\begin{equation*}
 	P_0 \in \Psi^0_0\lrpar{X;\bbE}
	\text{such that $[{}^0\sigma\lrpar{P_0}] = p.$}
\end{equation*}

 \item There is a fully elliptic element of the scattering calculus 
\begin{equation*}
 	P_{sc} \in \Psi^0_{\mathrm{sc}}\lrpar{X;\bbE}
	\text{such that $[{}^{sc}\sigma\lrpar{P_{sc}}] = p.$}
\end{equation*}

\item The Atiyah-Bott obstruction of $[p]$ vanishes, i\@.e\@.,
\begin{equation*}
p \in \ker\lrpar{ \Kc\lrpar{T^*X} \xrightarrow{r^*} \Kc\lrpar{T^*_{\pa X}X}}.
\end{equation*}
\end{enumerate}

Furthermore, each such realization determines a lift of the class $p$ to an
element of $\Kc\lrpar{T^*X,T_{\pa X}X},$ and the index is given by the
Atiyah-Singer index theorem applied to the double of $X.$
\end{theorem}

\begin{proof} The equivalence of $(1)$, $(2)$, and $(3)$ with $(4)$ has
already been established: for $(1)$ this is Theorem 5.14 in \cite{Boutet},
for $(2)$ this follows from \cite[(6.12)]{Melrose-Rochon}, and for $(3)$
this was established in Proposition \ref{TopSixTerm} above.  For each of
these elements of $\Kc\lrpar{T^*X}$ there are in fact many fully elliptic
elements of each calculus with this fixed symbol. This is
evident for the scattering calculus, and follows for Boutet de Monvel's
calculus and the zero calculus from the existence of a homotopy to an
operator whose symbol is the identity near the boundary (the homotopy of
the symbol restricted to the boundary cosphere bundle defines an extension
of the symbol).

The last part of the theorem follows from a version of Proposition 2.1 in \cite{Atiyah}. Set
\begin{multline*}
Q\lrpar{X,\pa X}=\\
\{ \lrpar{\alpha, \beta} \in \mathrm{ISO}_{S^*X}\lrpar{E,F} \times
\mathrm{ISO}_{B^*\pa X}\lrpar{E,F} : \alpha\rest{S^*\pa X} =
\beta\rest{S^*\pa X} \} 
\end{multline*}
with distinguished subset
\begin{equation*}
Q^*\lrpar{X}
=\{ \alpha \in \mathrm{ISO}_{S^*X}\lrpar{E,F} : \alpha\rest{S^*\pa X} = \Id \}.
\end{equation*}

\begin{lemma}
There exists a unique function
\begin{equation*}
	Q\lrpar{X,\pa X} \xrightarrow{f} K^0\lrpar{T^*X,T^*_{\pa X}X}
\end{equation*}
satisfying the following conditions:
\begin{itemize}
	\item [i)]
$\displaystyle f\lrpar{\alpha,\beta} = f\lrpar{\alpha',\beta'}$
if $\lrpar{\alpha,\beta}$ and $\lrpar{\alpha',\beta'}$ are homotopic in
$Q\lrpar{X,\pa X}.$ 
\item [ii)]
$\displaystyle f\lrpar{\lrpar{\alpha,\beta}\oplus\lrpar{\alpha',\beta'}} = 
f\lrpar{\alpha,\beta}+f\lrpar{\alpha',\beta'}$
\item [iii)]
$\displaystyle f\lrpar{\alpha,\Id} = [\alpha]$ if $\alpha \in Q^*\lrpar{X}$.
\end{itemize}
\end{lemma}

\begin{proof} This follows directly from the analysis of the scattering
case in \cite{Melrose-Rochon}.  Indeed, stable homotopy classes of elements
in $Q\lrpar{X,\pa X}$ are precisely the classes in $\curly
K^0_{\mathrm{sc}}$, so the existence follows from the isomorphism
\begin{equation*}
\curly K^0_{\mathrm{sc}} \cong K^0\lrpar{T^*X,T^*_{\pa X}X}.
\end{equation*}
 
Uniqueness follows from an argument of Atiyah. Any $\lrpar{\alpha,
  \beta} \in Q\lrpar{X,\pa X}$ is stably homotopic to an element in
  $Q^*\lrpar{X},$ \ie, there is an element $\lrpar{\alpha',\Id}\in
  Q^*\lrpar{X}$ with $f\lrpar{\alpha',\Id}=0$ and a homotopy of 
  $\lrpar{\alpha,\beta} \oplus \lrpar{\alpha',\Id}$ to another
  $\lrpar{\alpha'',\Id} \in Q^*\lrpar{X},$ hence
\begin{equation*}
	f\lrpar{\alpha,\beta} = 
	f\lrpar{\lrpar{\alpha,\beta}\oplus\lrpar{\alpha',\Id}}
	= f\lrpar{\alpha'',\Id} = [\alpha''].
\end{equation*}
\end{proof}
\end{proof}

\subsection{Differential operators}

Given a local elliptic boundary value problem of Lopatinskii-Schapiro
type, $\lrpar{P,B},$ it follows from the discussion above that there is a
Fredholm zero operator with the same interior symbol and index. In fact, if
$P$ is a differential operator of order $k$, then $x^kP$ is an elliptic
zero-differential operator whose interior symbol coincides with that of $P$. It is also possible to give an
explicit example of a perturbation of order $-\infty$ to a fully elliptic
operator corresponding to $B.$

\begin{theorem} Let $(P,B)$ be a local, fully elliptic boundary value
  problem, then there exists an operator $\curly{B}\in
  \Psi_0^{-\infty}\lrpar{X;E,F}$ such that $x^kP + \curly{B} \in
  \Psi^k_0\lrpar{X;E,F}$ is a fully elliptic operator with the same index as 
$\lrpar{P,B}$.
\end{theorem}

\begin{proof} By assumption, $P$ is an elliptic differential operator between
sections of $E$ and $F$ and $B$ is a formal differential operator at the
boundary between sections of $E$ and $G.$ The normal symbol of $P$ can be
used to identify $E$ and $F$ near the boundary, so we can assume that
$E\sim F$ near $\pa X$ and that the normal symbol of $P$ is the identity (for
some choice of normal vector field). Then let $p\lrpar{x\pa_x,x\hat\eta}$ and
$b(x\pa_x,x\hat\eta)$ be the model operators obtained from the principal
symbols frozen at the boundary. We shall show that the Lopatinskii-Schapiro
condition is equivalent to the statement that if $\chi(x)$ is a cutoff near
$x=0$ then for $\eps>0$ small enough
\begin{equation}
\langle b^*g,\chi ^2\lrpar{\frac x\eps }bf\rangle
\label{6.4.2006.1}\end{equation}
is non-degenerate as a pairing between the null space of $p$ and the null
space of $p^*.$ From this it follows by continuity that the operator of
order $-\infty$ in the b-calculus on $[0,\infty)$
\begin{equation}
a=I_\delta b^*\chi ^2(\frac x\eps )b I_\delta
\label{6.4.2006.2}\end{equation}
is such that $p+a$ is invertible for small $\delta;$ here $I_\delta$ is an
approximate identity. Then $a$ is in the range of the reduced normal map
and $x^mP+\eps A$ is Fredholm for small $\eps >0$ if
\begin{equation}
\sigma (A)=a.
\label{6.4.2006.3}\end{equation}

To see the non-degeneracy of \eqref{6.4.2006.1} observe that
the Lopatinskii-Schapiro condition is the condition that
\begin{equation*}
	u \xrightarrow{\beta^+} b\lrpar{\pa_x,\hat\eta}u\big|_{x=0}
\end{equation*}
is an isomorphism from the $L^2$ null space of $P(\pa_x,\hat\eta)$ to the
fibre of the bundle $G$ at the corresponding point. This in turn is
equivalent to the non-degeneracy of the corresponding bilinear form on the
product of the null space of $p$ and of $p^*$
\begin{equation}
  \langle b^*v,bu\rangle _G.
\label{6.4.2006.4}\end{equation}
The null space of $p\lrpar{\pa_x,\hat\eta }$ is the same as that of
$p\lrpar{x\pa_x,x\hat\eta}$ and the difference between \eqref{6.4.2006.1} and
\eqref{6.4.2006.4} is that there is an integral in \eqref{6.4.2006.1}, cut
off by $\chi$ and with additional factors of $x.$ It follows that
\eqref{6.4.2006.1} has an asymptotic expansion as $\eps \downarrow0$
with leading term (with power determined by the orders, \ie, the powers of
$x)$ with coefficient just \eqref{6.4.2006.4}. Thus, for small $\eps$
the form \eqref{6.4.2006.1} is indeed non-degenerate. Inserting the
approximate identities as in \eqref{6.4.2006.2} gives a bilinear form
converging to the previous one as $\delta \to0$ so again this is
non-degenerate, now for $\delta >0$ small enough.

The resulting operator $a$ is then in the range of the reduced normal
operator, e.g., its b-symbol is independent of $\hat\eta$.

It follows that $\pi'a\pi$ is an isomorphism where $\pi$ and $\pi'$ are
the projections onto the null space of $p$ and that of $p^*.$ An element in
the null space of the operator $p+\eps a$ can be decomposed according
to these projections, into $u+v$ with $u$ in the null space of $p$ and $v$
orthogonal to it. Then it must satisfy
\begin{equation}
(p+\epsilon a_{11})u_1+\epsilon a_{12}v=0,\
a_{21}u+a_{22}v=0
\label{6.4.2006.5}\end{equation}
where $a_{ij}$ is the $2\times2$ decomposition of $a$ with respect to these
projections. By arrangement, $a_{22}$ is invertible so this reduces to
\begin{equation}
\left(p+\eps (a_{11}-a_{12}a_{22}^{-1}a_{21})\right)u=0.
\label{6.4.2006.6}\end{equation}
For $\eps >0$ this has no solutions, by the invertibility of $p$ as a
map from the orthocomplement of its null space to its range. Thus the
reduced normal operator is invertible and hence the operator is $P+\epsilon
A$ Fredholm.

Finally, recall from \cite{Atiyah} that the isomorphism $\beta^+,$ hence the
quadratic form \eqref{6.4.2006.4}, determines the index and indeed the lift
of the K-class of $\sigma\lrpar{P}$ from $K\lrpar{T^*X}$ to
$K\lrpar{T^*X,T^*\pa X}$. Thus the index of $x^kP + \curly{B}$ is the same
as that of $\lrpar{P,B}$.
\end{proof}

\def\cprime{$'$} \def\cprime{$'$} \def\cdprime{$''$} \def\cprime{$'$}
  \def\cprime{$'$} \def\cprime{$'$} \def\cprime{$'$} \def\bud{$''$}
  \def\cprime{$'$} \def\cprime{$'$} \def\cprime{$'$} \def\cprime{$'$}
  \def\cprime{$'$} \def\cprime{$'$} \def\cprime{$'$} \def\cprime{$'$}
  \def\cprime{$'$} \def\cprime{$'$}
  \def\polhk#1{\setbox0=\hbox{#1}{\ooalign{\hidewidth
  \lower1.5ex\hbox{`}\hidewidth\crcr\unhbox0}}} \def\cprime{$'$}
  \def\cprime{$'$} \def\cprime{$'$} \def\cprime{$'$} \def\cprime{$'$}
  \def\cprime{$'$} \def\cprime{$'$} \def\cprime{$'$} \def\cprime{$'$}
  \def\cprime{$'$} \def\cprime{$'$} \def\cprime{$'$} \def\cprime{$'$}
  \def\cprime{$'$} \def\cprime{$'$} \def\cprime{$'$} \def\cprime{$'$}
  \def\cprime{$'$} \def\cprime{$'$} \def\cprime{$'$} \def\cprime{$'$}
  \def\cprime{$'$} \def\cprime{$'$} \def\cprime{$'$} \def\cprime{$'$}
  \def\cprime{$'$} \def\cprime{$'$}
\providecommand{\bysame}{\leavevmode\hbox to3em{\hrulefill}\thinspace}
\providecommand{\MR}{\relax\ifhmode\unskip\space\fi MR }
\providecommand{\MRhref}[2]{%
  \href{http://www.ams.org/mathscinet-getitem?mr=#1}{#2}
}
\providecommand{\href}[2]{#2}


\end{document}